\documentclass[12pt]{article}
\usepackage{graphicx, mathrsfs,amssymb,amsmath,amsfonts,amsbsy,amsthm,latexsym,bm,float}
\begin{document}


\title{
A fast algorithm for the linear canonical transform
\vskip.5cm}
\author{Rafael G. Campos and Jared Figueroa\\
Facultad de Ciencias F\'{\i}sico-Matem\'aticas,\\
Universidad Michoacana, \\
58060, Morelia, Mich., M\'exico.\\
\hbox{\small rcampos@umich.mx, jared@fismat.umich.mx}\\
}
\date{}
\maketitle
{
\noindent Keywords: Linear Canonical Transform, Fractional Fourier Transform, Quadrature, Hermite polynomials, Fractional Discrete Fourier Transform,  {\tt fft}
}\\
\vspace*{1truecm}
\begin{center} Abstract \end{center}
In recent years there has been a renewed interest in finding fast algorithms to compute accurately the linear canonical transform (LCT) of a given function. This is driven by the large number of applications of the LCT in optics and signal processing. The well-known integral transforms: Fourier, fractional Fourier, bilateral Laplace and Fresnel transforms are special cases of the LCT. 
In this paper we obtain an ${\mathcal O}(N\log N)$ algorithm to compute the LCT by using a chirp-FFT-chirp transformation yielded by a convergent quadrature formula for the fractional Fourier transform. This formula gives a unitary discrete LCT in closed form. In the case of the fractional Fourier transform the algorithm computes this transform for arbitrary complex values inside the unitary circle and not only at the boundary. In the case of the ordinary Fourier transform the algorithm improves the output of the FFT.
\vskip1.5cm
\section{Introduction}\label{intro}
The Linear Canonical Transform (LCT) of a given function $f(x)$ is a three-parameter integral transform that was obtained in connection with canonical transformations in Quantum Mechanics \cite{Mos71, Wol79}. It is defined by
\begin{eqnarray}\label{lctuno}
{\mathcal L}^{\{a,b,c,d\}}[f(x),y]&=&\frac{1}{\sqrt{2\pi i  b}}\int_{-\infty}^\infty e^{\frac{i}{2b}(ax^2-2xy+dy^2)}f(x)dx,\nonumber
\end{eqnarray}
for $b\ne0$, and by $\sqrt{d} e^{\frac{i}{2}cdy^2} f(dy)$, if $b=0$. The four parameters $a$, $b$, $c$ and $d$ appearing in (\ref{lctuno}), are the elements of a $2\times 2$ matrix with unit determinant, i.e., $ad-bc=1$. Therefore, only three parameters are free. Since this transform is a useful tool for signal processing and optical analysis, its study and direct computation in digital computers has become an important issue \cite{Hea09}-\cite{Pei05}, particularly, fast algorithms to compute the linear canonical transform have been devised \cite{Koc08, Hen05}. These algorithms use the following related ideas: (a) use of the periodicity and shifting properties of the discrete LCT to break down the original matrix into smaller matrices as the FFT does with the DFT, (b) decomposition of the LCT into a chirp-FFT-scaling transformation and (c) decomposition of the LCT into a fractional Fourier transform followed by a scaling-chirp multiplication. All of these are algorithms of ${\mathcal O}(N\log N)$ complexity.\\
In this paper we present an algorithm that takes ${\mathcal O}(N\log N)$ time based in the decomposition of the LCT into a scaling-chirp-DFT-chirp-scaling transformation, obtained by using a quadrature formula of the continuous Fourier transform \cite{Cam09, Cam92}.  Here, DFT stands for the standard discrete Fourier transform. To distinguish this discretization from other implementations, we call it  the extended Fourier Transform (XFT). 
Thus, the quadrature from which the XFT is obtained, uses some asymptotic properties of the Hermite polynomials and yields a fast algorithm to compute the Fourier transform, the fractional Fourier transform and therefore, the LCT. The quadrature formula is ${\mathcal O}(1/N)$-convergent to the continuous Fourier transform for certain class of functions \cite{Cam95}.
\section{A discrete fractional Fourier transform}\label{secdos}
In previous work \cite{Cam92}, \cite{Cam95}, \cite{Cam08},  we derived a quadrature formula for the continuous Fourier transform which yields an accurate discrete Fourier transform. For the sake of completeness we give in this section a brief review of the main steps to obtain this formula.\\
Let us consider the family of Hermite polynomials $H_n(x)$, $n=0,1,\ldots$, which satisfies the recurrence equation 
\begin{equation}\label{receqg}
H_{n+1}(x)+2nH_{n-1}(x)=2xH_n(x),
\end{equation}
with $H_{-1}(x)\equiv 0$. 
Note that the recurrence equation (\ref{receqg}) can be written as the eigenvalue problem 
\begin{equation}\label{eighinf}
\begin{pmatrix} 0&1/2&0&\cdots \\1&0& 1/2&\cdots \\0& 2& 0&\cdots \\\vdots&\vdots&\vdots&\ddots\\ \end{pmatrix}
\begin{pmatrix} H_0(x)\\  H_1(x)\\  H_2(x)\\ \vdots\end{pmatrix}=x\begin{pmatrix} H_0(x)\\  H_1(x)\\  H_2(x)\\ \vdots\end{pmatrix}.
\end{equation}
Let us now consider the eigenproblem associated to the principal submatrix of dimension $N$ of (\ref{eighinf})
\[
{\mathcal H}=\begin{pmatrix}0&1/2&0&\cdots & 0& 0\\1& 0& 1/2&\cdots & 0& 0\\0& 2& 0&\cdots &  0& 0\\
\vdots&\vdots&\vdots&\ddots&\vdots&\vdots\\
0& 0& 0&\cdots & 0&1/2\\0& 0& 0&\cdots & N-1&0\end{pmatrix}.
\]
It is convenient to symmetrize ${\mathcal H}$ by using the similarity transformation $S{\mathcal H} S^{-1}$ where $S$ is the diagonal matrix 
\[
S=\text{diag}\left\{1,\frac{1}{\sqrt{2}},\ldots,\frac{1}{\sqrt{(N-1)!\,2^{N-1}}}\right\}.
\]
Thus, the symmetric matrix $H=S{\mathcal H} S^{-1}$ takes the form
\[
\begin{pmatrix}0&\sqrt{\frac{1}{2}}&0&\cdots & 0& 0\\\sqrt{\frac{1}{2}}& 0& \sqrt{\frac{2}{2}}&\cdots & 0& 0\\
0& \sqrt{\frac{2}{2}}& 0&\cdots &  0& 0\\
\vdots&\vdots&\vdots&\ddots&\vdots&\vdots\\
0& 0& 0&\cdots & 0&\sqrt{\frac{N-1}{2}}\\0& 0& 0&\cdots & \sqrt{\frac{N-1}{2}}&0\end{pmatrix}.
\]
The recurrence equation (\ref{receqg}) and the Christoffel-Darboux formula \cite{Sze75} can be used to solve the eigenproblem 
\[
Hu_k=x_k u_k,\quad k=1,2,\ldots, N,
\]
which is a finite-dimensional version of (\ref{eighinf}). The eigenvalues $x_k$ are the zeros of $H_N(x)$ and the $k$th eigenvector $u_k$ is given by 
\[
c_k\left(s_1 H_0(x_k),s_2 H_1(x_k),\cdots,s_N H_{N-1}(x_k)\right)^T,
\]
where $s_1,\ldots,s_N$ are the diagonal elements of $S$ and $c_k$ is a normalization constant that can be determined from the condition $u_k^T\,u_k=1$, i.e., from
\[
c_k^2\,\sum_{n=0}^{N-1} \frac{H_n(x_k)H_n(x_k)}{2^n n!}=1.
\]
Therefore,
\[
c_k=\sqrt{\frac{2^{N-1}\,(N-1)!}{N}}\,\frac{(-1)^{N+k}}{ H_{N-1}(x_k)}.
\]
Thus, the components of the orthonormal vectors $u_k$, $k=1,2,\ldots, N$, are
\begin {equation}\label{ortvec}
(u_k)_n=(-1)^{N+k}\sqrt{\frac{2^{N-n}\,(N-1)!}{N\,(n-1)!}}\,\frac{H_{n-1}(x_k)}{H_{N-1}(x_k)},
\end{equation}
$n=1,\ldots,N$. Let $U$ be the orthogonal matrix whose $k$th column is $u_k$ and let us define the matrix  
\[
{\mathcal F}_z=\sqrt{2\pi}U^{-1}D(z)U,
\]
where $D(z)$ is the diagonal matrix $D(z)=\text{diag}\{1,z,z^2,\ldots,z^{N-1}\}$ and $z$ is an complex number. Therefore, the components of ${\mathcal F}_z$ are given by
\begin{eqnarray}\label{tmat}
({\mathcal F}_z)_{jk}&=&\sqrt{2\pi}\,\frac{(-1)^{j+k}\,2^{N-1}\,(N-1)!}{N\,H_{N-1}(x_j)H_{N-1}(x_k)}\sum_{n=0}^{N-1}\frac{z^n}{2^n\,n!}H_n(x_j)H_n(x_k).
\end{eqnarray}
Next, we want to prove that if $N$ is large enough, (\ref{tmat}) approaches the kernel of the fractional Fourier transform evaluated at $x=x_j$, $y=x_k$. To this, we use the asymptotic expression for $H_N(x)$ \cite{Sze75})
\begin{equation}\label{asymhnt}
H_N(x)\simeq\frac{\Gamma(N+1)e^{x^2/2}}{\Gamma(N/2+1)}\cos(\sqrt{2N+1}\,\,x-\frac{N\pi}{2}).
\end{equation}
Thus, the asymptotic form of the zeros of $H_N(x)$ are
\begin{equation}\label{asymcer}
x_k= \left(\frac{2 k-N-1}{\sqrt{2N}}\right)\frac{\pi}{2},
\end{equation}
$k=1,2,\ldots,N$. The use of (\ref{asymhnt}) and (\ref{asymcer}) yields
\[
H_{N-1}(x_k)\simeq (-1)^{N+k}\,\frac{\Gamma(N)}{\Gamma(\frac{N+1}{2})}\,e^{x_k^2/2},\quad N\to\infty,
\]
and the substitution of this asymptotic expression in (\ref{tmat}) yields
\begin{eqnarray*}
({\mathcal F}_z)_{jk}&\simeq&\sqrt{2\pi}\,\frac{2^{N-1}\,[\Gamma(\frac{N+1}{2})]^2}{\Gamma(N+1)}\,\,e^{-(x_j^2+x_k^2)/2}\sum_{n=0}^{\infty}\frac{z^n}{2^n\,n!}H_n(x_j)H_n(x_k).
\end{eqnarray*}
Finally, Stirling's formula and Mehler's formula \cite{Erd53} produce 
\begin{equation}\label{fasy}
({\mathcal F}_z)_{jk}\simeq \sqrt{\frac{2}{1-z^2}}\,e^{-\frac{(1+z^2)(x_j^2+x_k^2)-4x_jx_kz}{2(1-z^2)}}\Delta x_k,
\end{equation}
where $\Delta x_k$ is the difference between two consecutive asymptotic Hermite zeros, i.e.,
\begin{equation}\label{deltk}
\Delta x_k=x_{k+1}-x_k=\frac{\pi}{\sqrt{2N}}.
\end{equation}
Let us consider now the vector of samples of a given function $f(x)$
\[
f=(f(x_1),f(x_2),\ldots,f(x_N))^T.
\] 
The multiplication of the matrix ${\mathcal F}_z$ by the vector $f$ gives the vector $g$ with entries
\begin{eqnarray}\label{primcuad}
g_j&=&\sum_{k=1}^N({\mathcal F}_z)_{jk}f(x_k)\simeq \sqrt{\frac{2}{1-z^2}}\sum_{k=1}^N e^{-\frac{(1+z^2)(x_j^2+x_k^2)-4x_jx_kz}{2(1-z^2)}}f(x_k)\Delta x_k,\nonumber
\end{eqnarray}
where $j=1,2\ldots,N.$ This equation is a Riemann sum for the integral
\begin{eqnarray}\label{eqifft}
{\mathscr F}_z[f(x),y]&=&\sqrt{\frac{2}{1-z^2}}\int_{-\infty}^\infty e^{-\frac{(1+z^2)(y^2+x^2)-4xyz}{2(1-z^2)}}f(x)dx,\nonumber 
\end{eqnarray}
where $\vert z\vert<1$. Therefore, if we make $y_j=x_j$,
\begin{equation}\label{cuadtrfor}
{\mathscr F}_z[f(x),y_j]\simeq\sum_{k=1}^N({\mathcal F}_z)_{jk}f(x_k),\quad N\to\infty.
\end{equation}
Note that ${\mathscr F}_z[g(x),y]$ is the continuous fractional Fourier transform \cite{Nam80} of $g(x)$ except for a constant and therefore, ${\mathcal F}_z$ is a discrete fractional Fourier transform.  
\section{A fast linear canonical transform}\label{sectres}
Firstly, note that if $b\ne 0$, the LCT can be written as a chirp-FT-chirp transform
\begin{eqnarray}\label{lctdos}
{\mathcal L}^{\{a,b,c,d\}}[f(x),y]&=&\frac{e^{\frac{idy^2}{2b}}}{\sqrt{2\pi i  b}}\int_{-\infty}^\infty e^{-\frac{ixy}{b}}e^{\frac{iax^2}{2b}}f(x)dx.\nonumber
\end{eqnarray}
Thus, for $b\ne 0$, the LCT of the function $f(x)$ can be represented by the $1/b$-scaled Fourier transform of the function $g(x)=e^{\frac{iax^2}{2b}}f(x)$,
multiplied by $\frac{e^{\frac{idy^2}{2b}}}{\sqrt{2\pi i  b}}$.\\
On the other hand, note that for the case $z=\pm i$, (\ref{fasy}) yields a discrete Fourier transform $({\mathcal F}_{\pm i})_{jk}\simeq e^{\pm it_jt_k}\Delta t_k$, that can be related to the standard DFT as follows. The use of (\ref{asymcer}) yields
\begin{equation}\label{fftuno}
({\mathcal F}_i)_{jk}=e^{\pm it_jt_k}\Delta t_k=\frac{\pi}{\sqrt{2N}} e^{i\frac{\pi^2}{2N} \left(j-\frac{N-1}{2}\right) \left(k-\frac{N-1}{2}\right)}
\end{equation}
where we have used (\ref{asymcer}) and (\ref{deltk}). Since $\sum_{k=1}^N({\mathcal F}_i)_{jk}f(x_k)$ is a quadrature and therefore, an approximation of
\[
g(y_j)=\int_{-\infty}^\infty e^{i y_j x} f(x)dx,
\]
a scaled Fourier transform
\begin{equation}\label{furesc}
\int_{-\infty}^\infty e^{i \kappa y_j x} f(x)dt=g(\kappa y_j)
\end{equation}
has the quadrature $\sum_{k=1}^N\tilde{F}_{jk}f(x_k)$, where 
\begin{equation}\label{fftdos}
\tilde{F}_{jk}=\frac{\pi}{\sqrt{2N}}e^{i \kappa\frac{\pi^2}{2N} \left(j-\frac{N-1}{2}\right) \left(k-\frac{N-1}{2}\right)}.
\end{equation}
If we choose $\kappa=4/\pi$, (\ref{fftdos}) takes the form 
\begin{eqnarray}\label{dislctn}
F_{jk}&=&\frac{\pi e^{i\frac{\pi}{2}\frac{(N-1)^2}{N}}}{\sqrt{2N}}\left[e^{-i\pi\frac{N-1}{N} j}\right]\left[e^{i\frac{2\pi}{N}jk}\right]\nonumber
\left[e^{-i\pi\frac{N-1}{N} k}\right],
\end{eqnarray}
for $ j,k=0,1,2,\ldots,N-1$, and $\sum_{k=1}^NF_{jk}f(x_k)$ is an approximation of $g(4 y_j/\pi)$.  If now we choose $\kappa=4b/\pi$, but we keep the same matrix (\ref{dislctn}), then $\sum_{k=1}^NF_{jk}f(x_k)$ is an approximation of 
\[
\int_{-\infty}^\infty e^{i\frac{y_j}{b} x} f(x)dt.
\]
If now we replace $f(x)$ by $e^{\frac{iax^2}{2b}}f(x)$ and take into account (\ref{lctdos}), we have that 
\[
\sum_{k=1}^NF_{jk}e^{iax_k^2/2b}f(x_k)
\]
is an approximation of the product of functions $ \left(\frac{e^{\frac{idy^2}{2b}}}{\sqrt{2\pi i  b}}\right)^{-1}{\mathcal L}^{\{a,b,c,d\}}[f(x),y]$ evaluated at $y_j=4bx_j/\pi$.  Therefore, a discrete (scaled) linear canonical transform $L$ can be given in closed form. If we denote by $G(y)$ the LCT of $f(x)$, then
\[
G(y_j)=G(4bx_j/\pi)=\sum_{k=1}^N(S_1FS_2)_{jk}f(x_k),
\]
where $S_1$ and $S_2$ are diagonal matrices whose diagonal elements are $e^{\frac{idy_j^2}{2b}}/\sqrt{2\pi i  b}$, and $e^{iax_j^2/2b}$, respectively. As it can be seen, the matrix $L=S_1FS_2$, which gives the discrete LTC, i.e., the XFT, consists in a chirp-DFT-chirp transformation, where DFT stands for the standard discrete Fourier transform. Therefore, we can use any FFT to give a fast computation of the linear canonical transform $G=Lf$.\\
Now, the fast algorithm for the linear canonical transform can be given straightforwardly. 

\fbox{
\begin{minipage}{15cm}
\begin{center} \vskip 0.4cm
Algorithm \end{center}
\hrule width 15cm \vskip.3truecm
To compute an approximation $G=(G_1,G_2,\ldots,G_N)^T$ of the linear canonical transform $G(y)={\mathcal L}^{\{a,b,c,d\}}[f(x),y]$ evaluated at $y_j=4bx_j/\pi$, where $x_j= \left(\frac{2 j-N-1}{\sqrt{2N}}\right)\frac{\pi}{2}$.

\begin{enumerate}
\item For given $N$ set up the vector $u$ of components 
\begin{eqnarray*}
u_k&=&e^{-i\pi \frac{(k-1)(N-1)}{N}}e^{iax_k^2/2b}f\left(\pi \frac{2 k-N-1}{2\sqrt{2N}}\right),
\end{eqnarray*}
$k=1,2,\ldots,N$.
\item Set $y_j=4bx_j/\pi$ and compute the diagonal matrix $S$ according to 
\[
S_{jk}=\frac{\pi e^{i\frac{\pi}{2}\frac{(N-1)^2}{N}}}{\sqrt{2N}}\frac{e^{\frac{idy_j^2}{2b}}e^{-i\pi\frac{N-1}{N}(j-1)}}{\sqrt{2\pi i  b}}\delta_{jk},
\]
$j,k=1,2,\ldots,N$.
\item Let $D_F$ be the discrete Fourier transform, i.e., $(D_F)_{jk}=e^{i\frac{2\pi}{N}jk}$, $j,k=0,1,2,\ldots,N-1$. Obtain the approximation $G_j$ to $G(\frac{4b}{\pi} x_j)$ by computing the matrix-vector product
\begin{equation}\label{algfxft1}
G=SD_Fu,
\end{equation}
with a standard FFT algorithm.
\end{enumerate}
\end{minipage}}
\section{Example}
For this example we take an integral formula found in \cite{Gra94} that gives
\begin{eqnarray}\label{exa1}
G(y)&=&\frac{\sqrt{\pi}}{\sqrt{2\pi ib}(\alpha^2+\frac{a^2}{4 b^2})^{1/4}} e^{\frac{\alpha(\beta^2-\alpha\gamma)}{\alpha^2+\frac{a^2}{4 b^2}}}
e^{\frac{i}{2} \arctan(\frac{a}{2\alpha b})} e^{-\frac{\alpha y^2+2 \beta a y+a^2\gamma}{4 b^2\alpha^2+a^2}} \\
&\times& e^{i\frac{a c y^2}{2(4 b^2 \alpha^2+a^2)}} e^{i\frac{2 b (\alpha^2 d y^2+2 \beta\alpha y+\beta^2 a)}{4 b^2 \alpha^2+a^2}},\nonumber
\end{eqnarray}
if $f(x)=e^{-(\alpha y^2+2 \beta y+\gamma)}$, $\alpha>0$.
Figure 1 shows the exact LTC with $\alpha=1$, $\beta=2$, $\gamma=3$, $a=1$, $b=2$, $c=1/2$, and $d=2$, compared with the approximation given by the XFT.
\begin{figure}\label{Figu}
\centering
\includegraphics[scale=0.6]{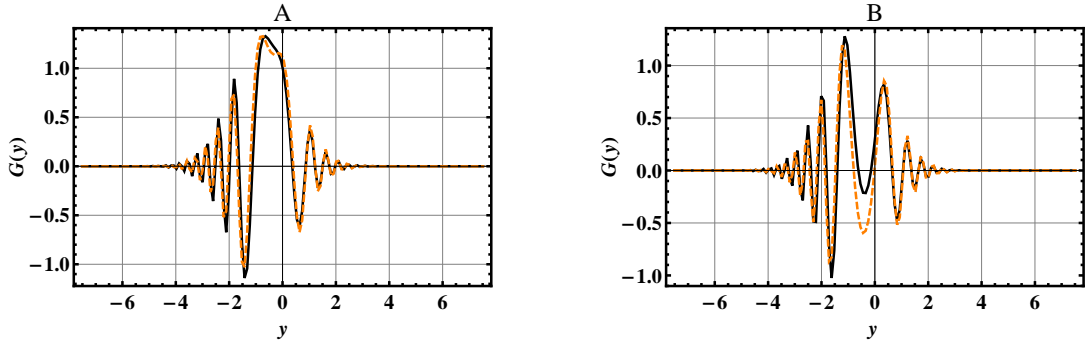}
\caption{Real part ({\bf A}) and imaginary part ({\bf B}) of the exact linear canonical transform (solid line) compared with the output of the XFT (dashed line) computed with $N=512$. The function $G(y)$ is that given in (\ref{exa1}) for $\alpha=1$, $\beta=2$, $\gamma=3$, $a=1$, $b=2$, $c=1/2$, and $d=2$.}
\end{figure}
Figure 2 shows the exact Fresnel transform of $f(x)=e^{-(\alpha y^2+2\beta  y+\gamma)}$ with $\alpha=2$, $\beta=1$, $\gamma=3$, $a=1$, $b=100$, $c=0$, and $d=1$, compared with the approximation given by the XFT.
\begin{figure}\label{Figd}
\centering
\includegraphics[scale=0.55]{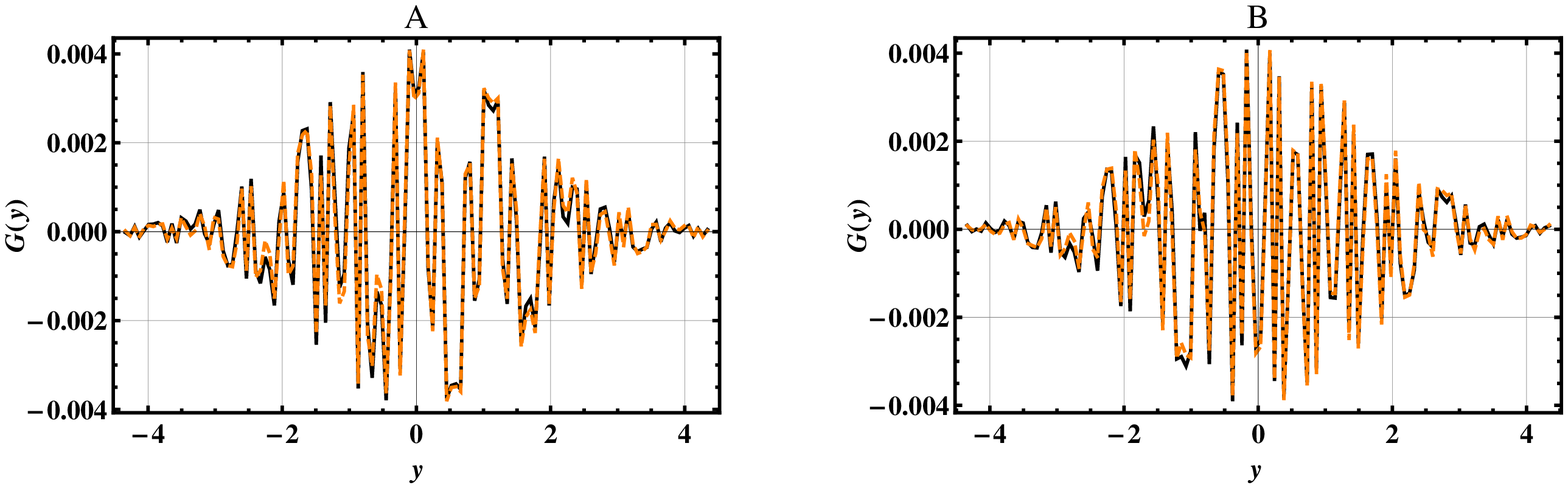}
\caption{Real part ({\bf A}) and imaginary part ({\bf B}) of the exact Fresnel transform (solid line) compared with the output of the XFT (dashed line) computed with $N=1024$. The function $G(y)$ is that given in (\ref{exa1}) for $\alpha=2$, $\beta=1$, $\gamma=3$, $a=1$, $b=100$, $c=0$, and $d=1$.}
\end{figure}
\section{Conclusion}\label{seccuatro}
We have obtained a discrete linear canonical transform and a fast algorithm to compute this transform by projecting the space of functions onto a vector space spanned by a finite number of Hermite functions. The XFT is a discrete LCT given by a unitary matrix in a closed form in which the DFT can be found at the core, surrounded by diagonal transformations, which makes easy to implement it in a fast algorithm. Since this discrete LCT is related to a  quadrature formula of the fractional Fourier transform, it yields accurate results.


\end{document}